\documentclass[jsco]{academic}
\usepackage{amsmath,amssymb,bbm}
\usepackage[light]{draftcopy}
\usepackage[arrow,curve,matrix]{xy}

\def\C{\mathcal{C}}
\def\ellv{\mathbf{l}}
\def\ev{\mathbf{e}}
\def\fv{\mathbf{f}}
\def\G{\mathcal{G}}
\def\I{\mathcal{J}}
\def\iv{\mathbf{i}}
\def\jv{\mathbf{j}}
\def\kv{\mathbf{k}}
\def\kk{{\mathbbm k}}
\def\L{\mathcal{L}}
\def\lcm#1{\mathrm{lcm}(#1)}
\def\lc#1#2{\mathrm{lc}_{#1}#2}
\def\lspan#1{\langle{#1}\rangle}
\def\lt#1#2{\mathrm{lt}_{#1}#2}
\def\M{\mathcal{M}}
\def\one{{\mathbbm 1}}
\def\P{\mathcal{P}}
\def\R{\mathcal{R}}

\def\S{\mathcal{S}}
\def\Sv{\mathbf{S}}
\def\syz#1#2{\mathrm{Syz}^{#1}(#2)}
\def\Spoly#1#2#3{S_{#1}(#2,#3)}
\def\Spolym#1#2#3{\Sv_{#1}(#2,#3)}
\def\T{\mathcal{T}}
\def\V{\mathcal{V}}

\def\sdr#1#2#3#4#5{%
  \xymatrix{{#1}\ar@<2pt>@/^7pt/[rr]^{#3} && {#2}
    \ar@<2pt>@/^7pt/[ll]^{#4} \ar@(dr,ur)[]_{#5}}}

\begin{document}
\shorttitle{Taylor and Lyubeznik Resolutions via Gr\"obner Bases}
\shortauthor{W.M. Seiler}
\title{Taylor and Lyubeznik Resolutions via Gr\"obner Bases}
\author{Werner M. Seiler}
\address{Lehrstuhl f\"ur Mathematik I, Universit\"at Mannheim\\ 
     68131 Mannheim, Germany\\
     Email: \texttt{werner.seiler@math.uni-mannheim.de}\\
    Web: 
    \texttt{www.math.uni-mannheim.de/\raisebox{-3pt}{\symbol{126}}wms}}
\maketitle

\begin{abstract}
  Taylor presented an explicit resolution for arbitrary monomial ideals.
  Later, Lyubeznik found that already a subcomplex defines a resolution. We
  show that the Taylor resolution may be obtained by repeated application of
  the Schreyer Theorem from the theory of Gr\"obner bases, whereas the
  Lyubeznik resolution is a consequence of Buchberger's chain criterion.
  Finally, we relate Fr\"oberg's contracting homotopy for the Taylor complex
  to normal forms with respect to our Gr\"obner bases and use it to derive a
  splitting homotopy that leads to the Lyubeznik complex.
\end{abstract}

\section{The Taylor and the Lyubeznik Resolution}\label{sec:taylor}

Let $\M=\{m_1,\dots,m_r\}\subset\P=\kk[x_1,\dots,x_n]$ be a finite set of
monomials. \citet{dt:res} constructed in her Ph.D. thesis an explicit free
resolution of the monomial ideal $\I=\lspan{\M}$. The associated complex
consists essentially of an exterior algebra and a differential defined via the
least common multiples of subsets of~$\M$.

Let $\V$ be some $r$-dimensional $\kk$-vector space with the basis
$\{v_1,\dots,v_r\}$.  If $\kv=(k_1,\dots,k_q)$ is a sequence of integers with
$1\leq k_1<k_2<\cdots<k_q\leq r$, we set
$m_{\kv}=\lcm{m_{k_1},\dots,m_{k_q}}$. The $\P$-module
$\T_q=\P\otimes\Lambda^q\V$ is then freely generated by all wedge products
$v_{\kv}=v_{k_1}\wedge\cdots\wedge v_{k_q}$.  Finally, we introduce on the
algebra $\T=\P\otimes\Lambda\V$ the following $\P$-linear differential
$\delta$:
\begin{equation}\label{eq:detaylor}
  \delta v_{\kv}=\sum_{\ell=1}^q(-1)^{\ell-1}
                     \frac{m_{\kv}}{m_{\kv_\ell}}v_{\kv_\ell}\;,
\end{equation}
where $\kv_\ell$ denotes the sequence $\kv$ with the entry $k_\ell$ removed.
Obviously, the differential $\delta$ respects the grading of $\T$ by the form
degree, as it maps the component $\T_q$ into $\T_{q-1}$ (however, in general
$\delta$ does not respect the natural \emph{bi}grading of $\T$ given by
$\T_{rq}=\P_r\otimes\Lambda^q\V$).

One can show that $(\T,\delta)$ is a complex representing a free resolution of
the ideal $\I_0=\lspan{\M}$.  Obviously, $\delta v_i=m_i$ and the length of
the resolution is given by the number $r$ of monomials. This implies
immediately that the resolution is rarely minimal.\footnote{Several
  characterisations of the case that $\T$ defines a minimal resolution have
  been given by \citet{rf:complex}.} Note that the ordering of the monomials
$m_i$ in the set $\M$ has no real influence on the result: the arising
resolutions are trivially isomorphic.

\citet{gl:res} proved later in his Ph.D. thesis that in fact already a
subcomplex $\L\subseteq\T$ defines a free resolution of $\I_0$. Let $\kv$
again be an integer sequence; we denote for $1\leq i<r$ by $\kv_{>i}$ the
subsequence of all entries $k_j>i$. If we eliminate from the basis of the
Taylor complex $\T$ all generators $v_{\kv}$ where for at least one $1\leq
i<r$ the monomial $m_i$ divides $m_{\kv_{>i}}$, then the remaining part $\L$
is still a complex defining a resolution of $\I$.

Here the ordering of the monomials $m_i$ is crucial; in general, for different
orderings different eliminations will be possible. As the Taylor complex is
essentially independent of the orderings, one also obtains a free
subresolution of $\T$ via a ``reverse'' form of the Lyubeznik approach.
Namely, we define $\kv_{<i}$ as the subsequence of all entries $k_j<i$ and
then eliminate all generators $v_{\kv}$ where for at least one $1\leq i<r$ the
monomial $m_i$ divides $m_{\kv_{<i}}$.

Both resolutions are of considerable interest in homological algebra
\citep{jls:polyres,wms:spencer}.  \citet{rf:complex} constructed an explicit
\emph{contracting homotopy} for the Taylor complex that also restricts to the
Lyubeznik subcomplex, i.\,e.\ a $\kk$-linear map
$\psi:\T_q\rightarrow\T_{q+1}$ such that $\delta\psi+\psi\delta=1$. Given a
term $x^\mu v_{\kv}\in T_q$, let $\iota=\iota(x^\mu v_{\kv})$ be the minimal
value for $i$ such that $m_i\mid x^\mu m_{\kv}$. Then we define
\begin{equation}\label{eq:psi}
  \psi(x^\mu v_{\kv})=[\iota<k_1]
      \frac{x^\mu m_{\kv}}{m_{(\iota,\kv)}}v_{(\iota,\kv)}
\end{equation}
where $(\iota,\kv)$ denotes the sequence $(\iota,k_1,\dots,k_q)$ and $[\cdot]$
is the Kronecker-Iverson symbol \citep{gkp:concrete} which is $1$, if the
contained condition is true and $0$ otherwise.

\section{Gr\"obner Bases and Syzygies}

Gr\"obner bases \citep{al:gb,bw:groe,clo:iva} are an important tool in
computational algebra and have been introduced in the Ph.D. thesis of
\citet{bb:diss}.  If $\I$ is an ideal in the polynomial ring
$\P=\kk[x_1,\dots,x_n]$, then a finite set $\G\subset\I$ is a \emph{Gr\"obner
  basis} of the ideal $\I$ for a term order $\prec$, if the leading term
$\lt{\prec}{f}$ of any polynomial $f\in\I$ is divisible by the leading term
$\lt{\prec}{g}$ of a generator $g\in\G$. If we write
$\lt{\prec}{\I}=\lspan{\{\lt{\prec}{f}\mid f\in\I\}}$ for the monomial ideal
generated by the leading terms of all the elements of $\I$, then we may
express this defining condition concisely as
$\lspan{\lt{\prec}{\G}}=\lt{\prec}{\I}$.

If (and only if) the set $\G$ is a Gr\"obner basis of the ideal $\I$ for the
term order~$\prec$, then every polynomial $f\in\I$ possesses a so-called
\emph{standard representation} $f=\sum_{g\in\G}P_gg$ where all the polynomials
$P_g\in\P$ satisfy $\lt{\prec}{(P_gg)}\preceq\lt{\prec}{f}$. Of course, this
representation is not unique, as one may add arbitrary syzygies. More
generally, every polynomial $f\in\P$ may be written (with the help of the
so-called \emph{division algorithm}) in the form
\begin{equation}\label{eq:nf}
  f=\sum_{g\in\G}P_gg+\hat f
\end{equation}
where no term of $\hat f$ is contained in $\lt{\prec}{\I}$. The polynomial
$\hat f$ is called the \emph{normal form} of $f$ with respect to $\G$. One can
show that $\hat f$ is uniquely defined, if and only if $\G$ is a Gr\"obner
basis.

A central concept in the theory of Gr\"obner bases is that of an
\emph{$S$-polynomial} (the $S$ stands for syzygy). Given two polynomials
$f_1,f_2\in\P$, we define their $S$-polynomial $\Spoly{\prec}{f_1}{f_2}$ as
follows. Let $m_i=\lt{\prec}{f_i}$ be the leading term of $f_i$ and set
$m_{12}=\lcm{m_1,m_2}$. Then
\begin{equation}
  \Spoly{\prec}{f_1}{f_2}=
      \frac{m_{12}}{\lc{\prec}{(f_1)}m_1}f_1-
      \frac{m_{12}}{\lc{\prec}{(f_2)}m_2}f_2
\end{equation}
where $\lc{\prec}{(f_i)}$ denotes the leading coefficient of $f_i$.  An
important criterion for a set~$\G$ to be a Gr\"obner basis of the ideal
$\lspan{\G}$ is that for all $g_1,g_2\in\G$ the normal form of the
$S$-polynomial $\Spoly{\prec}{g_1}{g_2}$ vanishes.

All these notions generalise trivially to submodules of free $\P$-modules. We
consider elements of a free $\P$-module of rank $m$ as $m$-dimensional vectors
with polynomial entries. The $S$-polynomial is then defined to be zero, if the
two leading terms $\lt{\prec}{\fv_1}$ and $\lt{\prec}{\fv_2}$ belong to
different components.

\citet{fos:syz} proved in his diploma thesis that the standard representations
of the $S$-polynomials lead not only to a generating set of the first syzygy
module but in fact again to a Gr\"obner basis with respect to a special term
order $\prec_{\G}$ induced by the basis $\G=\{g_1,\dots,g_r\}$: let
$\{\ev_1,\dots,\ev_r\}$ denote the standard basis of a free $\P$-module whose
rank is $|\G|=r$; then we define for arbitrary monomials $s,t\in\P$ that
$s\ev_\alpha\prec_{\G}t\ev_\beta$, if either
$\lt{\prec}{(sg_\alpha)}\prec\lt{\prec}{(tg_\beta)}$ or both
$\lt{\prec}{(sg_\alpha)}=\lt{\prec}{(tg_\beta)}$ and $\beta<\alpha$.

Let $\G=\{g_1,\dots,g_r\}$ and $m_i=\lt{\prec}{g_i}$. Assume that the standard
representation of the $S$-polynomial $\Spoly{\prec}{g_i}{g_j}$ is given by
$\Spoly{\prec}{g_i}{g_j}=\sum_{k=1}^rP_{ijk}g_k$.  Then we set again
$m_{ij}=\lcm{m_i,m_j}$ and define the syzygy
\begin{equation}
  \Sv_{ij}=\frac{m_{ij}}{\lc{\prec}{(g_i)}m_i}\ev_i-
           \frac{m_{ij}}{\lc{\prec}{(g_j)}m_j}\ev_j-
           \sum_{k=1}^rP_{ijk}\ev_k\in\P^r\,.
\end{equation}
By the \emph{Schreyer theorem}, the set $\Sigma_{\G}=\{\Sv_{ij}\mid1\leq
i<j\leq r\}$ is a Gr\"obner basis of the syzygy module $\syz{}{\G}$ with
respect to the term order $\prec_{\G}$. We also introduce the syzygy
$\tilde{\Sv}_{ij}=\frac{m_{ij}}{m_i}\ev_i-\frac{m_{ij}}{m_j}\ev_j$ of the
leading terms of $g_i$ and $g_j$. As any monomial set is trivially a Gr\"obner
basis of the ideal it generates, we conclude that the set
$\tilde\Sigma_{\G}=\{\tilde{\Sv}_{ij}\mid1\leq i<j\leq r\}$ is a Gr\"obner
basis of the syzygy module $\syz{}{\lt{\prec}{\G}}$ with respect to the term
order $\prec_{\G}$.

\emph{Buchberger's chain criterion} \citep{bb:crit} asserts that certain
$S$-poly\-no\-mials may be ignored in the above mentioned criterion for a
Gr\"obner basis. It is based on the following observation which will turn out
to be crucial for the Lyubeznik resolution.  Let $\S\subseteq\Sigma_\G$ and
assume that (i) the set $\tilde\S=\{\tilde{\Sv}_{ij}\mid\Sv_{ij}\in\S\}$
generates the syzygy module $\syz{}{\lt{\prec}{\G}}$, (ii) $\tilde\S$ contains
the three syzygies $\tilde{\Sv}_{ij}$, $\tilde{\Sv}_{ik}$, $\tilde{\Sv}_{jk}$,
and (ii) the monomial $m_i$ divides $m_{jk}$. Then the set
$\tilde\S\setminus\{\tilde{\Sv}_{jk}\}$ still generates
$\syz{}{\lt{\prec}{\G}}$ and the set $\S\setminus\{\Sv_{jk}\}$ still generates
$\syz{}{\G}$.

\section{Taylor Resolution via Schreyer Theorem}

Our goal is now to show that the Taylor resolution can be constructed via
repeated application of the Schreyer theorem. The decisive point will be to
define an appropriate ordering of the generators in each Gr\"obner basis.

As a first step, we introduce on the free polynomial module $\T_q$ defined in
Sect.~\ref{sec:taylor} two term orders $\prec_q$ and $\prec_q^r$ as follows.
We define on the space of all ascending integer sequences a ``lexicographic''
order: we set $\kv<\ellv$, if for $j=\min\{i\mid k_i\neq\ell_i\}$ the
inequality $k_j<\ell_j$ holds.  For $\prec_0=\prec_0^r$ we choose an arbitrary
term order on $\P$; it will turn out that everything we do in the sequel is
independent of this choice.  Then we define recursively for two sequences
$\kv$ and $\ellv$ of length $q+1$ and two monomials $s,t\in\P$ that
$su_{\kv}\prec_{q+1}tu_{\ellv}$, if either $\lt{\prec_{q}}{(s\delta
  u_{\kv})}\prec_{q} \lt{\prec_{q}}{(t\delta u_{\ellv})}$ or both
$\lt{\prec_{q}}{(s\delta u_{\kv})}= \lt{\prec_{q}}{(t\delta u_{\ellv})}$ and
$\kv<\ellv$. For the ``reverse'' term order $\prec_{q+1}^r$ the last condition
is replaced by $\kv>\ellv$

\begin{lemma}\label{lem:lttaylor}
  If\/ $\kv$ is a sequence of length\/ $q$ with\/ $1\leq q\leq r$, then 
  \begin{equation}\label{eq:prec}
    \frac{m_{\kv}}{m_{\kv_q}}v_{\kv_q}\prec_{q-1}\cdots\prec_{q-1}
    \frac{m_{\kv}}{m_{\kv_1}}v_{\kv_1}
  \end{equation}
  and hence the leading term of\/ $\delta v_{\kv}$ is given by
  \begin{equation}\label{eq:ltdv}
    \lt{\prec_{q-1}}{(\delta v_{\kv})}=\frac{m_{\kv}}{m_{\kv_1}}v_{\kv_1}\;.
  \end{equation}
  For the reverse order $\prec_{q-1}^r$ we obtain\/
  $\tfrac{m_{\kv}}{m_{\kv_1}}v_{\kv_1}\prec_{q-1}^r\cdots\prec_{q-1}^r
  \tfrac{m_{\kv}}{m_{\kv_q}}v_{\kv_q}$ and thus\/ $\lt{\prec_{q-1}^r}{(\delta
    v_{\kv})}=\frac{m_{\kv}}{m_{\kv_q}}v_{\kv_q}$.
\end{lemma}

\begin{proof}
  We proceed by induction. For $q=1$ the assertion is trivial. For $q=2$ we
  must compare the two terms
  $\tfrac{m_{\kv}}{m_{\kv_1}}v_{\kv_1}=\tfrac{m_{\kv}}{m_{k_2}}v_{k_2}$ and
  $\tfrac{m_{\kv}}{m_{\kv_2}}v_{\kv_2}=\tfrac{m_{\kv}}{m_{k_1}}v_{k_1}$. As
  $\delta v_i=m_i$, we find that $\tfrac{m_{\kv}}{m_{k_1}}\delta
  v_{k_1}=\tfrac{m_{\kv}}{m_{k_2}}\delta v_{k_2}=m_{\kv}$. By the definition
  of the term order $\prec_1$ this implies that
  $\tfrac{m_{\kv}}{m_{k_2}}v_{k_2}\prec_1\tfrac{m_{\kv}}{m_{k_1}}v_{k_1}$ as
  $k_1<k_2$ (independent of the choice of the order $\prec_0$ as claimed
  above).
  
  Now assume that the lemma holds for all sequences of length less than $q$.
  We must compare the terms $\tfrac{m_{\kv}}{m_{\kv_\ell}}v_{\kv_\ell}$ with
  $1\leq\ell\leq q$. This requires to determine the leading term of
  $\tfrac{m_{\kv}}{m_{\kv_\ell}}\delta v_{\kv_\ell}$ with respect to
  $\prec_{q-2}$. By our induction hypothesis, we obtain for $\ell=1$ the term
  $\tfrac{m_{\kv}}{m_{\kv_{1,2}}}v_{\kv_{1,2}}$ and for $\ell>1$ the term
  $\tfrac{m_{\kv}}{m_{\kv_{1,\ell}}}v_{\kv_{1,\ell}}$. Thus the values
  $\ell=1$ and $\ell=2$ yield both the same term and, as obviously
  $\kv_2<\kv_1$, we find $\tfrac{m_{\kv}}{m_{\kv_2}}v_{\kv_2}\prec_{q-1}
  \tfrac{m_{\kv}}{m_{\kv_1}}v_{\kv_1}$.
  
  In order to compare the terms for the other possible values of $\ell$, we
  must descend recursively. At the next lower form degree we find that
  $\ell=2$ and $\ell=3$ yield the same term, so that by the definition of our
  term order, the one for $\ell=2$ is the greater one. Continuing until degree
  one we obtain (\ref{eq:prec}).

  The proof for the reverse term order $\prec_{q-1}^r$ proceeds completely
  analogously. This time at form degree $q'$ the terms for $\ell=q'$ and
  $\ell=q'-1$ coincide and by definition of the order the one for $\ell=q'$ is
  greater.
\end{proof}

Now that we know the leading terms of the elements $\delta v_{\kv}$, the next
step is to show that with respect to both introduced term orders the set
\begin{equation}
  \Delta_q=\{\delta v_{\kv}\mid\kv=(k_1,\dots,k_{q+1})\}
\end{equation}
is a Gr\"obner basis of the submodule $\I_q$ it generates in $\T_q$.  Note
that we have $\I_q=\delta(\T_{q+1})$ because of the $\P$-linearity of
$\delta$.

\begin{lemma}\label{lem:syztaylor}
  Let\/ $\kv$ be an ascending integer sequence of length\/ $q$ with\/ $0\leq
  q<r$ and\/ $i,j$ two further integers with\/ $1\leq i<j<k_1$. Then the
  syzygy induced by the $S$-polynomial $\Spolym{\prec_q}{\delta
    v_{(i,\kv)}}{\delta v_{(j,\kv)}}$ is
  \begin{equation}\label{eq:syztaylor}
    \Sv_{(j,\kv),(i,\kv)}=\delta v_{(i,j,\kv)}\;.
  \end{equation}
  Similarly, if $i,j$ are two integers with\/ $k_q<i<j\leq r$, then the syzygy
  induced by the $S$-polynomial $\Spolym{\prec_q^r}{\delta v_{(\kv,i)}}{\delta
    v_{(\kv,j)}}$ is
  \begin{equation}\label{eq:syzrtaylor}
    \Sv_{(\kv,i),(\kv,j)}=(-1)^{q+1}\delta v_{(\kv,i,j)}\;.
  \end{equation}
\end{lemma}

\begin{proof}
  We prove the assertion only for $\prec_q$, as the proof for $\prec_q^r$
  proceeds again completely analogously. The sign in (\ref{eq:syzrtaylor})
  stems from the fact that there the sequence is manipulated at its end and
  the sign of the last summand in (\ref{eq:detaylor}) depends on the length of
  the sequence.
  
  In order to simplify the notation, we write $\iv=(i,\kv)$, $\jv=(j,\kv)$
  and, finally, $\bar\kv=(i,j,\kv)$. As $(\T,\delta)$ is a complex, obviously
  $\delta^2v_{\bar\kv}=0$ and, by the $\P$-linearity of $\delta$, we find that
  \begin{equation}\label{eq:dsquare}
    \sum_{\ell=1}^{q+2}(-1)^{\ell-1}
        \frac{m_{\bar\kv}}{m_{\bar\kv_\ell}}\delta v_{\bar\kv_\ell}=0\;.
  \end{equation}
  Thus $\delta v_{\bar\kv}\in\syz{}{\Delta_q}$. On the other hand, by
  definition of an $S$-polynomial and Lemma \ref{lem:lttaylor},
  $\Spolym{\prec_q}{\delta v_{\jv}}{\delta
    v_{\iv}}=\frac{m_{\bar\kv}}{m_{\jv}}\delta
  v_{\jv}-\frac{m_{\bar\kv}}{m_{\iv}}\delta v_{\iv}$. Comparing with
  (\ref{eq:dsquare}), we see that these are just the first two summands. Thus
  we are done, if we can show that the remaining terms of (\ref{eq:dsquare})
  define a standard representation of $\Spolym{\prec_q}{\delta v_{\jv}}{\delta
    v_{\iv}}$.
  
  $S$-polynomials are defined such that the leading terms of the two summands
  cancel in the subtraction. Hence in order to find the leading term of our
  $S$-polynomial we must compare the second largest term in each summand. A
  straightforward computation shows that applying $\delta$ yields in each case
  the same leading term. Hence we obtain
  \begin{equation}
    \lt{\prec_q}{\bigl(\Spolym{\prec_q}{\delta v_{\iv}}{\delta v_{\jv}}\bigr)}=
    \frac{m_{\bar\kv}}{m_{(j,\kv_1)}}v_{(j,\kv_1)}
  \end{equation}
  and there only remains to show that for all values $3\leq\ell\leq q+2$ the
  relation $\lt{\prec_q}{\bigl(\tfrac{m_{\bar\kv}}{m_{\bar\kv_\ell}} \delta
    v_{\bar\kv_\ell}\bigr)} \preceq_q
  \tfrac{m_{\bar\kv}}{m_{(j,\kv_1)}}v_{(j,\kv_1)}$ holds. This is
  straightforward, as Lemma \ref{lem:lttaylor} implies
  $\lt{\prec_q}{\bigl(\tfrac{m_{\bar\kv}}{m_{\bar\kv_\ell}} \delta
    v_{\bar\kv_\ell}\bigr)}=
  \tfrac{m_{\bar\kv}}{m_{(j,\kv_\ell)}}v_{(j,\kv_\ell)}$. Applying $\delta$
  and taking the leading term gives immediately the desired result. Now the
  lemma follows from the definition of the syzygy $\Sv_{\iv,\jv}$.
\end{proof}

\begin{proposition}\label{prop:gbtaylor}
  For every\/ $q$ with $0\leq q\leq r$, the set\/ $\Delta_q$ is a Gr\"obner
  basis both for\/ $\prec_q$ and for\/ $\prec_q^r$ of the submodule\/ $\I_q$.
\end{proposition}

\begin{proof}
  This proposition is a corollary to the above mentioned Schreyer theorem. We
  proceed again by induction. For $q=0$, the assertion is trivial, as
  $\Delta_0=\M$ and any monomial set is a Gr\"obner basis of the ideal
  generated by it for every term order.
  
  Let us assume that the proposition holds for $q-1$. In order to invoke the
  Schreyer theorem we must define an ordering of the elements $\delta
  v_{\kv}\in\Delta_q$. Two fairly natural choices are to order them either
  ascending or descending by the index $\kv$ (using the above defined
  lexicographic ordering of integer sequences).  It is easy to see that the
  term order $\prec_{\Delta_{q-1}}$ used in the Schreyer theorem is in the
  first case $\prec_q^r$ and in the second case $\prec_q$. The assertion is
  now a trivial consequence of Schreyer's theorem and the
  Lemmata~\ref{lem:lttaylor} and \ref{lem:syztaylor}.
\end{proof}

As an immediate corollary we arrive finally at the following theorem.

\begin{theorem}
  The Taylor resolution can be obtained by repeatedly applying the Schreyer
  theorem, if at each step the generators are ordered in either of the two
  ways described in the proof above.
\end{theorem}

\section{Lyubeznik Resolution and Buchberger Chain Criterion}

The Gr\"obner basis of the syzygy module obtained via the Schreyer theorem is
in general not reduced. A number of generators may be eliminated by applying
the Buchberger chain criterion. We show now that this allows us to derive the
Lyubeznik subcomplex $\L\subseteq\T$.

\begin{proposition}\label{prop:bbcrit}
  Assume that for some\/ $1\leq i<r$ and an integer sequence\/ $\kv$ of
  length\/ $q$ the monomial\/ $m_i$ divides\/ $m_{\kv_{>i}}$. Then the set\/
  $\Delta_q\setminus\{\delta v_{\kv}\}$ is still a Gr\"obner basis of the
  submodule\/ $\I_q$ for the term order\/ $\prec_q$.
\end{proposition}

\begin{proof}
  Suppose $m_i$ divides $m_{\kv_{>i}}$. We first assume that $\kv=\kv_{>i}$,
  i.\,e.\ $i<k_1$. The case $q=1$ is trivial: the Buchberger chain criterion
  implies that the syzygy $\Sv_{k_1,k_2}$ corresponding to the $S$-polynomial
  of $m_{k_1}$ and $m_{k_2}$ is not needed to generate the syzygy module
  $\syz{}{\M}$. By Lemma~\ref{lem:syztaylor}, $\Sv_{k_1,k_2}=\delta v_{\kv}$
  and thus $\Delta_1\setminus\{\delta v_{\kv}\}$ is still a Gr\"obner basis.
  
  For $q>2$ we consider the three elements $\delta v_{\kv_1}, \delta
  v_{\kv_2}, \delta v_{(i,\kv_{1,2})}\in\Delta_{q-1}$. According to Lemma
  \ref{lem:lttaylor}, all three leading terms are multiples of
  $v_{\kv_{1,2}}$; the coefficients are $\tfrac{m_{\kv_1}}{m_{\kv_{1,2}}}$,
  $\tfrac{m_{\kv_2}}{m_{\kv_{1,2}}}$ and
  $\tfrac{m_{(i,\kv_{1,2})}}{m_{\kv_{1,2}}}$, respectively. Our assumption
  $m_i\mid m_{\kv}$ implies that the least common multiple of the first two is
  divisible by the third. Thus we can invoke the Buchberger chain criterion
  and obtain our assertion, as $\Sv_{\kv_1,\kv_2}=\delta v_{\kv}$ according to
  Lemma \ref{lem:syztaylor}.
  
  If $\kv\neq\kv_{>i}$, then we may eliminate at an earlier stage the
  generator $\delta v_{\kv_{>i}}$. The generator $\delta v_{\kv}$ arises in
  our construction of the Taylor resolution by repeated application of the
  Schreyer theorem only as the result of chains of syzygies starting with one
  induced by an $S$-polynomial involving $\delta v_{\kv_{>i}}$. As the latter
  one may be eliminated, $\delta v_{\kv}$ is not needed either.
\end{proof}

\begin{proposition}
  Assume that for some\/ $1\leq i<r$ and an integer sequence\/ $\kv$ of
  length\/ $q$ the monomial\/ $m_i$ divides\/ $m_{\kv_{<i}}$. Then the set\/
  $\Delta_q\setminus\{\delta v_{\kv}\}$ is still a Gr\"obner basis of the
  submodule\/ $\I_q$ for the term order\/ $\prec_q^r$.
\end{proposition}

\begin{proof}
  Completely analogous to the previous proposition.
\end{proof}

\begin{theorem}
  The (reverse) Lyubeznik resolution arises from the Taylor resolution by
  repeated application of the Buchberger chain criterion.
\end{theorem}

Note that in general we do not invoke all possible instances of the Buchberger
chain criterion for deriving the Lyubeznik resolution. This follows already
from the trivial fact that the chain criterion is independent of the ordering
of the monomials in the original set $\M$ whereas the Lyubeznik criterion for
the elimination of generators of $\T$ is not.

\section{Contracting Homotopies and Normal Forms}

The contracting homotopy $\psi$ defined by (\ref{eq:psi}) also possesses an
interesting interpretation in terms of our Gr\"obner bases $\Delta_q$.

\begin{lemma}\label{lem:psi}
  We have $\psi(x^\mu v_{\kv})=0$, if and only if $x^\mu
  v_{\kv}\notin\lt{\prec_q}{\I_q}$. If, however, $x^\mu
  v_{\kv}\in\lt{\prec_q}{\I_q}$, then $\psi(x^\mu v_{\kv})=x^\nu
  v_{(\iota,\kv)}$ where $\iota$ is the minimal value of $i$ such that
  $\lt{\prec_q}{\delta v_{(i,\kv)}}\mid x^\mu v_{\kv}$ and $\nu$ is chosen
  such that $\lt{\prec_q}{(x^\nu\delta v_{(\iota,\kv)})}=x^\alpha v_{\kv}$.
\end{lemma}

\begin{proof}
  Let us first assume that in Fr\"oberg's definition $\iota(x^\mu v_{\kv})\geq
  k_1$ so that $\psi(x^\mu v_{\kv})=0$. According to Lemma~\ref{lem:lttaylor},
  the leading term of $\delta v_{\ellv}$ with respect to the term order
  $\prec_q$ only lies in the component generated by $v_{\kv}$, if
  $\ellv=(i,\kv)$ for some $1\leq i<k_1$. Thus we have $x^\mu
  v_{\kv}\in\lt{\prec_q}{\I_q}$, if and only if an exponent vector~$\nu$
  exists such that $x^\nu m_{(i,\kv)}=x^\alpha m_{\kv}$. But by the definition
  of the function $\iota$, the term $m_{(i,\kv)}$ does not divide $x^\alpha
  m_{\kv}$ for any $i<k_1$.
  
  Now let us assume that $\iota=\iota(x^\mu v_{\kv})<k_1$. Then it follows
  again from Lemma~\ref{lem:lttaylor} that $\lt{\prec_q}{\delta
    v_{(\iota,\kv)}}=\tfrac{m_{(\iota,\kv)}}{m_{\kv}}v_{\kv}$. By the
  definition of $\iota$, we have $m_{(\iota,\kv)}\mid x^\mu m_{\kv}$ and hence
  a unique exponent vector $\nu$ exists with $\lt{\prec_q}{(x^\nu\delta
    v_{(\iota,\kv)})}=x^\alpha v_{\kv}$. It is trivial that $\iota$ is the
  smallest value with this property.
\end{proof}

As a trivial corollary of this result we find that the restriction of the
contracting homotopy $\psi$ to the $\kk$-vector space generated by the terms
in $\lt{\prec_q}{\I_q}$ is injective. Probably of more interest is the
following observation.

\begin{theorem}\label{thm:psinf}
  The map\/ $\psi\circ\delta:\T_q\rightarrow\T_q$ yields the normal form with
  respect to the Gr\"obner basis\/ $\Delta_q$ and the map\/
  $\delta\circ\psi=1-\psi\circ\delta$ is a projector on\/ $\I_q$.
\end{theorem}

\begin{proof}
  Let $v\in\T_q$ be an arbitrary element. Applying the division algorithm with
  respect to the Gr\"obner basis $\Delta_q$ of $\I_q$ yields according to
  (\ref{eq:nf}) a representation
  $v=\sum_{|\kv|=q+1}\delta(P_{\kv}v_{\kv})+\hat v$ where $P_{\kv}\in\P$ are
  some polynomials and where $\hat v$ is the normal form of $v$ with respect
  to $\Delta_q$. Recall that the normal form is unique (in contrast to the
  coefficients $P_{\kv}$) and that it consists only of terms \emph{not}
  contained in $\lt{\prec_q}{\I_q}$.
  
  Exploiting Lemma~\ref{lem:psi} and that $\psi$ is a contracting homotopy, we
  find
  \begin{equation}
    \psi(v)=\sum_{|\kv|=q+1}\psi\delta(P_{\kv}v_{\kv})=
        -\sum_{|\kv|=q+1}\delta\psi(P_{\kv}v_{\kv})+
         \sum_{|\kv|=q+1}P_{\kv}v_{\kv}\;.
  \end{equation}
  This implies immediately
  $\delta\psi(v)=\sum_{|\kv|=q+1}\delta(P_{\kv}v_{\kv})$ and thus
  $\psi\delta(v)=\hat v$.
\end{proof}

In fact, we may extend this idea to a general principle for the construction
of contracting homotopies in complexes over the polynomial ring $\P$ where
Gr\"obner bases for the images of the differential are known.

\begin{theorem}\label{thm:conhom}
  Let\/ $(\C,\delta)$ be a (not necessarily finite) exact complex of free
  polynomial modules with a $\P$-linear differential\/
  $\delta:\C_q\rightarrow\C_{q-1}$.  Assume that for all degrees\/ $q$ a
  Gr\"obner basis of the submodule\/ $\delta(\C_q)$ is known.  Thus every
  element\/ $u\in\C_q$ may be written in the form\/ $u=\delta v+\hat u$
  where\/ $\hat u$ is the unique normal form of\/ $u$ with respect to the
  Gr\"obner basis of\/ $\delta(\C_q)$. Then the map\/
  $\psi:\C_q\rightarrow\C_{q+1}$ defined by\/ $\psi(u)=\hat v$ where\/ $\hat
  v$ is the unique normal form of\/ $v$ with respect to the Gr\"obner basis
  of\/ $\delta(\C_{q+1})$ is a contracting homotopy of\/ $(\C,\delta)$
  satisfying\/ $\psi^2=0$.
\end{theorem}

\begin{proof}
  We must first show that $\psi$ is well-defined. But this is trivial: if $v$
  and $v'$ are two elements of $\C_{q+1}$ such that $\delta v=\delta v'$, then
  $v=v'+\delta w$ by the exactness of $(\C,\delta)$ and thus $\hat v=\hat v'$.
  Then we must prove that $\psi\delta+\delta\psi=\one$. Again this is very
  simple, as, by definition of $\psi$, we find $\delta\psi(u)=\delta\hat
  v=\delta v$ and $\psi\delta(u)=\hat u$. Hence
  $(\psi\delta+\delta\psi)(u)=\delta v+\hat u=u$ as required. The relation
  $\psi^2=0$ follows trivially from the definition.
\end{proof}

It follows from Theorem~\ref{thm:psinf} that Fr\"oberg's contracting homotopy
is precisely the homotopy arising from this principle applied to the Gr\"obner
bases $\Delta_q$ with respect to the term orders $\prec_q$, although it
appears to be non-trivial (or at least rather tedious) to directly derive the
explicit expression (\ref{eq:psi}). 

We may also introduce a ``reverse'' contracting homotopy $\psi_r$ as follows.
For a given term $x^\mu v_{\kv}\in\T_q$ let $\iota=\iota(x^\mu v_{\kv})$ be
the maximal value for $i$ such that $m_i\mid x^\mu m_{\kv}$ and define
\begin{equation}
  \psi_r(x^\mu v_{\kv})=[\iota>k_q]
      \frac{x^\mu m_{\kv}}{m_{(\kv,\iota)}}v_{(\kv,\iota)}\;.
\end{equation}
This corresponds to applying Theorem~\ref{thm:conhom} to the Gr\"obner bases
$\Delta_q$ with respect to the ``reverse'' term orders $\prec_q^r$. We leave
the obvious details like the ``reverse'' form of Lemma~\ref{lem:psi} to the
reader.

Finally, we consider the restriction of the map $\psi$ to the Lyubeznik
complex~$\L$. From its definition (\ref{eq:psi}), it is not completely obvious
that $\psi(\L)\subset\L$ and thus that $\psi$ is a contracting homotopy for
$\L$, too; in fact, this was proven only very recently \citep{jls:polyres}.
Taking our approach, this becomes a simple corollary to Proposition
\ref{prop:bbcrit} and Theorem \ref{thm:conhom}.

\begin{corollary}
  The map\/ $\psi$ defined by (\ref{eq:psi}) is a contracting homotopy of the
  Lyubeznik complex\/ $\L$.
\end{corollary}

\begin{proof}
  The subset $\Delta_q^\prime=\{\delta(v_{\kv})\mid
  v_{\kv}\in\L_{q+1}\}\subseteq\Delta_q$ is, by Proposition \ref{prop:bbcrit},
  still a Gr\"obner basis of the ideal generated by $\Delta_q$. Thus applying
  Theorem \ref{thm:conhom} yields exactly the same contracting homotopy as for
  the full Taylor complex $\T$. This immediately implies that
  $\psi(\L)\subset\L$.
\end{proof}

\section{Strong Deformation Retracts}\label{sec:sdr}

\citet{bl:fpa} introduced the notion of a \emph{splitting homotopy} of a chain
complex $(\C,\delta)$ over a ring $\R$. This is a graded $\R$-module
homomorphism $\phi:\C\rightarrow\C$ such that (i)
$\phi(\C_q)\subseteq\C_{q+1}$ for all $q$, (ii) $\phi^2=0$, and (iii)
$\phi\delta\phi=\phi$. Such a map leads immediately to a \emph{strong
  deformation retract}
\begin{equation}\label{eq:sdr}
  \sdr{\S}{\C}{\iota}{\pi}{\phi}
\end{equation}
where $\pi=1-\delta\phi-\phi\delta$, $\S=\pi(\C)\subseteq\C$ and $\iota$ is
the inclusion map.  Indeed, it is easy to see that our assumptions imply
$\pi^2=\pi$ and one easily checks that $\pi\iota=\one_{\C}$ and
$\iota\pi=\pi$. One can show that any strong deformation retract arises
this way.

We are interested in the special case that $\I\subseteq\P$ is a polynomial
ideal and the chain complex $(\C,\delta)$ defines a free resolution of $\I$,
i.\,e.
\begin{equation}\label{eq:freeres}
  \cdots\longrightarrow\C_q\stackrel{\delta}{\longrightarrow}\C_{q-1}
  \longrightarrow\cdots\longrightarrow
  \C_0\stackrel{\epsilon}{\longrightarrow}
  \P/\I\longrightarrow 0
\end{equation}
is an exact sequence of free $\P$-module. Here we assume for simplicity that
$\C_0\subseteq\P$ and $\epsilon$ is the canonical projection. Obviously, both
the Taylor and the Lyubeznik resolution (or more generally any resolution
constructed via Schreyer's theorem) is of this form.  In such a situation,
strong deformation retracts allow us to construct smaller resolutions
\citep{lal:reshp}.

It is trivial to see that any contracting homotopy $\psi$ satisfying
$\psi^2=0$ is also a splitting homotopy. Thus, if the assumptions of Theorem
\ref{thm:conhom} are satisfied, then the map $\psi$ defined in it implies via
the isomorphism $\C_0/\delta_1(\C_1)\cong\P/\I$ a strong deformation
retract
\begin{equation}\label{eq:sdrpsi}
  \sdr{\P/\I}{\C}{c}{\epsilon}{\psi}
\end{equation}
where the projection\/ $\epsilon$ is extended to\/ $\C$ by setting it zero
outside of\/ $\C_0$ and where the map\/ $c$ yields the canonical
representative of each equivalence class with respect to the Gr\"obner basis
of\/ $\delta_1(\C_1)$.

However, the retract (\ref{eq:sdrpsi}) is not very interesting, as the left
hand side is not really a resolution anymore. We show now how the contracting
homotopy~$\psi$ of Theorem \ref{thm:conhom} may be used to obtain another
splitting homotopy $\phi$.  For this construction we generalise a homotopy
found by \citet{jls:polyres} (Proposition~3.9) which yields the Lyubeznik from
the Taylor resolution.

Let the vectors $\ev^{(q)}_\alpha$ with $1\leq\alpha\leq\dim\C_q$ form a basis of
the free module $\C_q$. Thus in the case of the Taylor resolution we may use
the vectors $v_{\kv}$ with $|\kv|=q$. Then we define inductively a map
$f:\C\rightarrow\C$ by setting $f(1)=1$,
\begin{equation}\label{eq:f}
  f(\ev^{(q)}_\alpha)=\psi\bigl(f(\delta\ev^{(q)}_\alpha)\bigr)
\end{equation}
and extending $\P$-linearly to the whole complex $\C$. The idea behind this
definition may be conveniently expressed in the following lemma.

\begin{lemma}
  Let\/ $\Delta_q$ be the Gr\"obner basis of\/ $\delta(\C_{q+1})$ and assume
  that we have\/ $\ev^{(q)}_\alpha\notin\lt{\prec}{\Delta_q}$ for all\/
  $1\leq\alpha\leq\dim\C_q$ and all\/ $q\geq0$. Then\/ $f$ is the identity map.
\end{lemma}

\begin{proof}
  Since $f(1)=1$, we find that
  $f(\ev^{(1)}_\alpha)=\psi\delta(\ev^{(1)}_\alpha)$. By definition of the
  contracting homotopy $\psi$, this means that $f(\ev^{(1)}_\alpha)$ is the
  normal form of $\ev^{(1)}_\alpha$ with respect to $\Delta_1$. If
  $\ev^{(1)}_\alpha\notin\lt{\prec}{\Delta_1}$, then this normal form is again
  $\ev^{(1)}_\alpha$. Now the assertion follows by induction.
\end{proof}

Note that $\ev^{(q)}_\alpha\in\lt{\prec}{\Delta_q}$ for some $\alpha$ and $q$
obviously implies that (\ref{eq:freeres}) is not the minimal resolution. The
map $f$ projects to a subcomplex $\S\subseteq\C$ by eliminating all generators
that appear in the leading ideals $\lt{\prec}{\Delta_q}$. Thus if $f$ is not
the identity map, we obtain a smaller resolution that is closer to the minimal
one.

Above we showed how a suitable projection $\pi:\C\rightarrow\S$ is obtained
from a splitting homotopy $\phi$. Now we have found a projection $f$ and would
like to get a suitable splitting homotopy $\phi$. As we know a contracting
homotopy for our complex, this is easily accomplished. We set $\phi(u)=0$
for\/ $u\in\C_0$. If we know $\phi$ restricted to $\C_{q-1}$, we can extend it
to $\C_q$ by solving the ``differential'' equation
$\delta\bigl(\phi(u)\bigr)=u-\phi\delta(u)-f(u)$. If we set
$v=u-\phi\delta(u)-f(u)$, then obviously $\delta(v)=0$ and we have
$\delta\psi(v)=v$. Thus a possible solution is $\phi(u)=\psi(v)$. This leads
to the following result.

\begin{theorem}\label{thm:splithom}
  Let the map\/ $\phi:\C\rightarrow\C$ be defined inductively by setting\/
  $\phi(u)=0$ for\/ $u\in\C_0$ and
  \begin{equation}\label{eq:phi}
    \phi(u)=\psi\bigl(u-\phi(\delta u)-f(u)\bigr)
  \end{equation}
  for\/ $u\in\C_q$ with\/ $q>0$. Then\/ $\phi$ is a splitting homotopy of the
  complex\/ $(\C,\delta)$ and we have a strong deformation retract
  \begin{equation}\label{eq:sdrphi}
    \sdr{\S}{\C}{\iota}{f}{\phi}
  \end{equation}
  where again\/ $\iota$ is the inclusion map.
\end{theorem}

\begin{proof}
  We have $\psi^2=0$ by Theorem \ref{thm:conhom} and this implies at once
  $\phi^2=0$. The relation $\phi\delta\phi=\phi$ follows by a simple
  computation from the fact that $\psi$ is a contracting homotopy and again
  that $\psi^2=0$.
\end{proof}

It is now easy to see that applying this construction to the Taylor complex
leads immediately to the Lyubeznik complex. Indeed, assume that $m_i\mid
m_{\kv}$ with $i<k_1$; then by Lemma \ref{lem:lttaylor} we find
\begin{equation}
  \lt{\prec_q}{\delta v_{(i,\kv)}}=
  \frac{m_{(i,\kv)}}{m_{\kv}}v_{\kv}=v_{\kv}\;.
\end{equation}
Thus the basis vector $v_{\kv}$ is not in image of $f$. But $m_i\mid m_{\kv}$
with $i<k_1$ is just Lyubeznik's condition for redundant generators in the
Taylor resolution. Obviously, the converse holds, too: if
$v_{\kv}\in\lt{\prec_q}{\Delta_q}$, then there exists an $i<k_1$ such that
$m_{(i,\kv)}=m_{\kv}$ which is equivalent to $m_i\mid m_{\kv}$.

\section{Conclusions}

We have demonstrated that both the Taylor and the Lyubeznik resolution are
directly derivable from the theory of Gr\"obner bases. The Taylor resolution
is a syzygy resolution obtained by applying the Schreyer theorem on a
Gr\"obner basis for the syzygy module. The Lyubeznik subresolution arises by
invoking some instances of the Buchberger chain criterion.

\citet{hmm:taylor,mm:hil,mm:res} studied in a series of articles the Taylor
resolution and discussed strategies to obtain smaller resolutions from it. One
of them lead to the Lyubeznik resolution.  Most (if not all) of their
reduction strategies may be homologically interpreted as splitting
homotopies.\footnote{This observation is due to Larry Lambe.}

The essential point in the proof of Proposition~\ref{prop:gbtaylor} that the
sets $\Delta_q$ are in fact Gr\"obner bases of the images of the differential
$\delta$ was to find the right ordering of the generators.  Actually, we found
two natural orderings on the bases $v_{\kv}$ that both lead to the Taylor
resolution.

The Schreyer theorem is often used to provide a simple proof of the Hilbert
syzygy theorem that every ideal in $\P=\kk[x_1,\dots,x_n]$ possesses a free
resolution of length at most $n$ (see e.\,g.\ \cite{al:gb}). Here yet another
ordering of the generators is used. Obviously, the length of the Taylor
resolution is in general much larger, as it is given by the number $s$ of
monomials in $\M$.  Thus in this respect our orderings based on the integer
sequences labelling the generators are not ``good''. This is the price to be
payed for the fact that the differential of the Taylor resolution has such a
simple explicit representation which makes it very useful in theoretical
considerations.

Theorems \ref{thm:conhom} and \ref{thm:splithom} demonstrate that the
contracting homotopy for the Taylor complex found by \citet{rf:complex} and
the splitting homotopy of \citet{jls:polyres} leading to the Lyubeznik
resolution, respectively, are not something particular but actually emerge
from general principles applicable to any complex of free polynomial modules.
Both theorems are not really surprising. One of the main tasks of a Gr\"obner
basis of an ideal $\I$ is to distinguish a unique representative in each
equivalence class in $\P/\I$ (this was the problem studied by
\citet{bb:diss}); a contracting homotopy does something fairly similar. Hence
a relation between the two concepts has to be expected.

The splitting homotopy constructed in Theorem \ref{thm:splithom} removes only
obstructions to the minimality of the resolution sitting in the leading terms
of the Gr\"obner basis of $\delta(\C_q)$. Another explicit resolution of
monomial ideals was recently derived with the help of Pommaret bases
\citep{wms:habil,wms:comb2}. While the construction yields automatically a
Gr\"obner basis of each image $\delta(\C_q)$, obstructions to minimality never
sit in their leading terms. Thus Theorem \ref{thm:splithom} yields only the
identity map in this case.

\section*{Acknowledgements}

The author is indebted to Larry Lambe for introducing him to the Taylor and
the Lyubeznik resolution. He also computed a number of useful examples with
the help of his \textsc{Axiom} programs for these resolutions. This work has
been financially supported by Deutsche Forschungsgemeinschaft and by {\small
  INTAS} (grant 99-1222).

\small
\bibliography{../../BIB/Algebra,../../BIB/Groebner,../../BIB/Misc,../../BIB/Seiler}
\bibliographystyle{plainnat}

\end{document}